\documentclass[reqno]{amsart}

\usepackage{amssymb,graphicx,color, amsthm}
\usepackage{stmaryrd}
\usepackage{atbegshi, picture}
\usepackage[colorlinks=true,citecolor=black,linkcolor=black,urlcolor=blue]{hyperref}
\usepackage[margin=1.2in, marginparwidth=2cm]{geometry}
\newcommand{\comment}[1]{}

\def\cVb{\cV_{\bullet}}
\def\cVw{\cV_{\circ}}
\def\Vb{V_{\bullet}}
\def\Vw{V_{\circ}}
\def\Esame{E^{=}}
\def\Edif{E^{\neq}}
\def\Ebb{E_{\bullet\bullet}}
\def\Ebw{E_{\bullet\circ}}
\def\Eww{E_{\circ\circ}}
\def\Ebbsame{\Ebb^=}
\def\Ebbdif{\Ebb^{\neq}}
\def\Ebwsame{\Ebw^=}
\def\Ebwdif{\Ebw^{\neq}}
\def\Ewwsame{\Eww^=}
\def\Ewwdif{\Eww^{\neq}}

\def\hPsi{\widehat{\Psi}}

\newcommand{\beqa}{\begin{eqnarray}}
\newcommand{\eeqa}{\end{eqnarray}}
\newcommand{\be}{\begin{equation}}
\newcommand{\ee}{\end{equation}}
\newcommand{\val}{\text{val}}

\def\cS{\mathcal{S}}

\theoremstyle{plain}
\newtheorem{theorem}{Theorem}[section]
\newtheorem{lemma}[theorem]{Lemma}
\newtheorem{coroll}[theorem]{Corollary}
\newtheorem{claim}[theorem]{Claim}
\newtheorem{prop}[theorem]{Proposition}

\theoremstyle{definition}
\newtheorem{definition}{Definition}[section]

\theoremstyle{remark}

\newcommand{\bGC}{\mathbb{G}_{\text{SYK}}}
\newcommand{\bGS}{\tilde{\mathbb{G}}_{\text{SYK}}}
\newcommand{\bS}{\mathbb{S}}
\newcommand{\bM}{\mathbb{M}}
\newcommand{\cG}{\mathcal{G}}
\newcommand{\cM}{\mathcal{M}}
\newcommand{\cR}{\mathcal{R}}
\newcommand{\cV}{\mathcal{V}}

\AtBeginShipoutNext{\AtBeginShipoutUpperLeft{
  \put(\dimexpr\paperwidth-0.5cm\relax,-0.8cm){\makebox[0pt][r]{\Small YITP-18-109}}
}}

\newcommand{\C}{c}

\def\cF{\mathcal{F}}

\newcommand{\boK}{\overline{\mathbb{K}}}

\newcommand{\cA}{\mathcal{A}}
\newcommand{\cT}{\mathcal{T}}
\newcommand{\cB}{\mathcal{B}}
\newcommand{\bG}{\mathbb{G}}

\newcommand{\bK}{\mathbb{K}}
\newcommand{\cK}{K}

\newcommand{\bC}{\mathbb{C}}
\newcommand{\bP}{\mathbb{P}}
\newcommand{\GF}{\mathcal{G}}
\newcommand{\Val}{\text{val}}

\def\GFbb{\GF_{\bullet\bullet}}
\def\GFbw{\GF_{\bullet\circ}}
\def\GFww{\GF_{\circ\circ}}

\def\zb{z_{\bullet}}
\def\zw{z_{\circ}}

\newcommand{\deltaS}{\delta}

\title[Combinatorial study of graphs arising from the SYK model]{Combinatorial study of graphs arising from the Sachdev-Ye-Kitaev model}

\author{\'E. Fusy, L. Lionni and A. Tanasa}

\date{\today}

\begin{document}

\begin{abstract} We consider the graphs involved in the theoretical physics model known as the colored Sachdev-Ye-Kitaev (SYK) model.
We study in detail their combinatorial properties  at any order in the so-called $1/N$ expansion, and we enumerate these graphs asymptotically.

Because of the duality between colored graphs involving $q+1$ colors and colored triangulations in dimension $q$, our results apply to the asymptotic enumeration of spaces that generalize unicellular maps - in the sense that they are obtained from a single building block - for which a higher-dimensional generalization of the genus is kept fixed. 

\end{abstract}

\maketitle


\section{Introduction}

In the last years, a quantum mechanical model known as the Sachdev-Ye-Kitaev (SYK) model \cite{SY, Kitaev} has attracted huge interest from the theoretical physics community (see \cite{PR,maldacena} and references within). This comes from the fact that the SYK model is the unique known model enjoying a certain number of important properties in an high-energy physics context, which makes it a pertinent toy-model for black hole physics. 

The SYK model (and related models) is studied ``perturbatively", that is, using formal decompositions of the quantities defining the theory in sums over graphs. 
In the case of the SYK model, one uses the  so-called $1/N$ expansion, where $N$ is the number of particles  of the model.
Some graphs give the \emph{leading} contribution (a first approximation in the limit $N\to\infty$), and the other graphs provide \emph{corrections} 
of relative importance to the leading contribution. This level of contribution of a graph (the importance of the correction it brings to the leading computation) is encoded in a non-negative parameter which we here call 
the \emph{order}. Until now, most papers focused on the leading contribution  \cite{PR, maldacena, nador}, whose corresponding graphs (the so-called melonic graphs) are very simple from a combinatorial point of view.

In this paper, we realize a combinatorial study of the graphs arising from the  perturbative $1/N$ expansion of the colored SYK model \cite{razvan, blt}, at {\it any order of contribution}. This is done using purely combinatorial methods: a bijection with combinatorial maps (constellations), and the method of kernel extraction \cite{Kernel1, Kernel2}.
The results of this paper are a follow up of preliminary studies initiated in \cite{blt, LL}.

In \cite{Witten}, Witten related the SYK model to the colored tensor model, a generalization of random matrix models initially proposed in a mathematical physics context by Gurau  \cite{Gurau,1N-arb} (see also the book \cite{Guraubook}), and then studied from a purely combinatorial point of view by Gurau and Schaeffer in \cite{gilles}\footnote{A similar combinatorial study, for a different type of tensor model, called the multi-orientable tensor model \cite{mo} was done in \cite{FT}.}.
Our analysis is analogous to the Gurau-Schaeffer study, even though the combinatorial objects we deal with here, called hereafter SYK graphs, are qualitatively different from the colored graphs analyzed in \cite{gilles}\footnote{More specifically, it is mainly the way these graphs are classified, their \emph{order}, which differs for both models.}. 
Because of these differences, the asymptotic analysis turns out to be significantly less involved for the colored SYK model than for the colored tensor model. 

\

Our work has another application in a discrete geometry context. Indeed, the graphs studied here are dual to colored triangulations \cite{Italians, GagliardiDipol, Chapuy-triang}, or gluings of bigger building blocks \cite{BLR, octa, LL}, and the \emph{order} mentioned  above is a certain linear combination of the number of sub-simplices of these spaces. In dimension two, and for surfaces obtained from a single polygon by identifying two-by-two the edges on its boundary (called unicellular maps \cite{unicell3, unicell1, unicell2}),
the order reduces to twice the genus of the surface. In higher dimension, we will show that asymptotically almost surely  (i.e.~with  probability tending to 1 when the number of vertices of the graphs goes to infinity; we use the abbreviation a.a.s.~thereafter), the spaces of fixed order are obtained from a single building block, whose boundary represents a piecewise-linear manifold, the topology of which we determine.

\subsection{Colored graphs and SYK colored graphs} 
\begin{figure}[!h]
\includegraphics[scale=0.8]{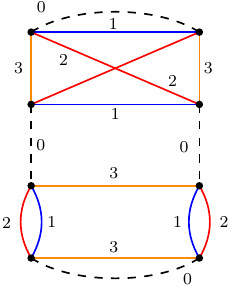} \hspace{1.5cm }\includegraphics[scale=0.9]{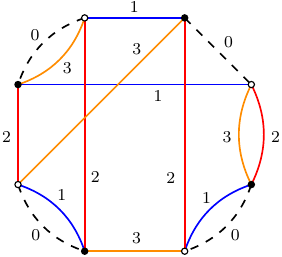}
\caption{\label{fig:ExamplesGraphs} A generic 4-colored graph and a bipartite SYK graph.}
\end{figure}

Let us start with the following definition:

\begin{definition}
\label{def:graph}
A connected regular $(q+1)$-edge-colored graph has edges carrying colors in $\{0,\cdots,q\}$, so that  each color reaches every vertex precisely once. Throughout the text, we will simply refer to such a graph as a {\it colored graph}. 

A colored graph is said to be {\it rooted} if one of its color-0 edges is distinguished and oriented. It is said to be {\it bipartite} if its vertices are colored in black and white so that every edge links a black and a white vertex. For a rooted bipartite graph we take the convention that the vertex at the origin of the root-edge is black. 
We denote by $\tilde \bG^q$ the family of connected rooted $(q+1)$-edge-colored graphs, and $\bG^q$ the subfamily of bipartite graphs from $\tilde \bG^q$. 

A connected colored graph is called an {\it SYK graph} if it remains connected when all the color-0 edges are deleted.
We denote by $\bGS^q$ the family of rooted $(q+1)$-edge-colored SYK graphs, and $\bGC^q$ the subset of bipartite SYK graphs. 
\end{definition}
We give in Fig.~\ref{fig:ExamplesGraphs} an example of a generic colored graph and of a bipartite SYK graph. 
As the color 0 will play a special role in the following, we represent the edges of color $0$ as dashed.

\subsection{Bi-colored cycles and order of a graph}
Throughout the paper, we will use the notation $\llbracket a,b\rrbracket = \{a,\ldots, b\}$ for two integers $a<b$. The valency of a vertex is the number of incident edges.

\begin{definition}
If $i\in\llbracket 1,q \rrbracket$, the connected components of the subgraph obtained by keeping only the edges of color 0 and $i$ only have vertices of valency two, so that they are cycles which alternate colors $0$ and $i$.
We will call such a subgraph a {\it color-$0i$ cycle}, 
and more generally a bi-colored cycle.\footnote{In the tensor model literature, these are referred to as ``faces". This terminology follows that of matrix models (where faces of ribbon graphs are closed cycles). However we choose not to use it to avoid confusions with the faces of maps and the facets of simplices.
} 
\end{definition}

For $i=1, \ldots, q$, we denote by $F_{0i}(G)$ the number of  color-$0i$ cycles of a colored graph $G$, and 
\begin{equation}
\label{eq:BicCyclesDef}
F_0(G) = \sum_{i=1}^q F_{0i}(G).
\end{equation}

\begin{definition}
Given a colored graph $G\in\bG^q$, we call {\it 0-residues} of $G$ the connected components of the graph $G_{\hat 0}$ obtained  from $G$ by deleting all the color-0 edges. We denote by $R_0(G)$ the number of its 0-residues.
\end{definition}

An SYK graph is therefore a colored graph $G$ with a single 0-residue, $R_0(G)=1$.

\begin{definition}
Denoting by $V(G)$ the number of vertices of a colored graph $G$, we define its {\it order},
\begin{equation}
\label{eq:DeltaS0}
\deltaS_0(G) = 1+ \frac{q-1} 2 V(G) - F_0(G).
\end{equation}
\end{definition}

It is known that this parameter is always non-negative (we will present a bijection to certain diagrams in Section~\ref{sec:Bijections}, which makes it easily visible).

\medskip 

\subsection{Feynman graphs of the SYK model}
\label{subsec:SYK}

The SYK graphs of Def.~\ref{def:graph} are the graphs that label the perturbative $1/N$ expansion of the colored SYK model of \cite{razvan, blt}\footnote{Which is a particular case of the Gross-Rosenhaus SYK model \cite{GR}.} (called Feynman graphs in physics).

More precisely, the 
so-called
(normalized) two-point function of the colored SYK model, one of the fundamental objects defining the theory, admits a formal $1/N$  expansion of the form 
\begin{equation}
\label{eq:ExpansionSYK}
\cG(N) =  \sum_{G\in\bGS^q}N^{ - \delta_0(G)}A(G).
\end{equation}
The sum above is taken over SYK graphs  (or bipartite SYK graphs for the complex SYK model), 
and $A(G)$ is a quantity which depends importantly on the details of the graph $G$.  We call this quantity the amplitude of $G$.
These amplitudes are only known for the very simplest cases  \cite{ PR, maldacena}.  The other fundamental objects defining the theory admit similar expansions.

Because $\deltaS_0$ is non-negative, this sum re-organizes as follows, 
\begin{equation}
\label{eq:ExpansionSYK2}
\cG(N) = \sum_{\delta\ge 0}\frac 1 {N^\delta }  \sum_{\substack{{G\in\bGS^q,}\\{\delta_0(G) = \delta}}}A(G 
). 
\end{equation}
The parameter $N$ is taken to be very large in a first approximation\footnote{Physically, it represents the number of particles (fermionic fields) described by the model.}. The graphs of order 0 provide the leading contribution, for which the amplitudes $A(G)$  can be computed \cite{SY, PR}, the graphs of order 1 provide the first correction in $1/N$, and so on:
the graphs of order $\deltaS_0$ provide the corrections of order $1/N^{\deltaS_0}$ to the large $N$ computations. 

In order to compute the amplitudes $A(G)\rvert_{\delta_0(G) = \delta}$, one must first identify the SYK graphs contributing at order $\deltaS_0(G) = \delta$. This was done for graphs of order zero \cite{PR} and one \cite{blt, DartoisSYK}. The procedure was described for higher orders in \cite{blt} and further detailed in \cite{LL}.

As already explained above, we are interested in the present paper in studying the combinatorial properties of the SYK graphs of any fixed order $\delta\geq 0$. 
We stress that because of the non-combinatorial amplitudes $A(G)$, which a priori differ for graphs contributing to the same order $\deltaS_0 = \delta$, our enumerative results do not allow re-summations of contributions to the SYK two-point functions. Our interest is therefore purely combinatorial.

\medskip

\subsection{Geometric interpretation: enumeration of unicellular discrete spaces}
\label{sub:Geom-Interp}

By duality,  edge-colored graphs $G$ with $q+1$ colors encode $q$-dimensional colored triangulations $\cT(G)$ (for more details on this duality, see e.g.~\cite{Italians}). The colored graph $G$ is said to represent the topological space $\lvert\cT(G)\rvert$  triangulated by $\cT(G)$, and every homeomorphic space.
A colored graph is bipartite if and only if the space it represents is orientable.  A $q$-dimensional  triangulation is colored if its $(q-1)$-simplices carry a  color $c\in \llbracket 0,q\rrbracket$, so that the $q+1$ $(q-1)$-simplices incident to a given $q$-simplex have different colors.  In our case however, the color 0 is given a special role, as we only focus on color-$0i$ cycles. In this case, we can interpret the colored graphs as dual to discrete spaces obtained by gluing together some building blocks along the elements of their triangulated $(q-1)$-dimensional boundaries \cite{LL}.
\begin{figure}[h!]
\includegraphics[scale=0.5]{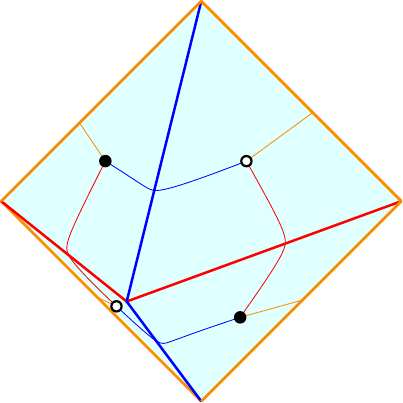}
\hspace{2cm}
\raisebox{15pt}{\includegraphics[scale=0.55]{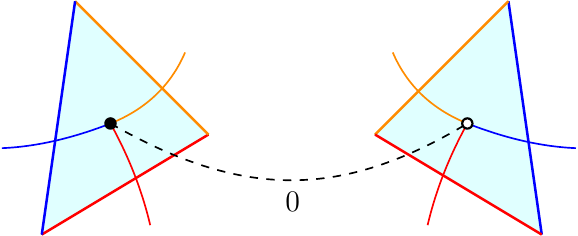}}
\caption{\label{fig:Unicell}On the left is an octahedron with the 3-colored graph dual to its boundary. 
To build a unicellular space out of the octahedron, one adds color-0 edges to identify two-by-two the 2-simplices of its boundary (triangles) in a unique way.}
\end{figure}

As a first example, studied in details in \cite{octa}, let us consider an octahedron in dimension $q=3$ (Fig.~\ref{fig:Unicell}). Its 2-dimensional boundary is triangulated, so that the edges (the colored 1-simplices of the triangulated boundary) carry colors in $\{1,2,3\}$, and so that the triangles are incident to edges of all colors in $\{1,2,3\}$. This 2-dimensional colored triangulation is dual to a 3-colored graph whose edges carry colors in $\{1,2,3\}$ (left of Fig.~\ref{fig:Unicell}). By adding edges of color 0 between the black and white vertices of this 3-colored graph (right of Fig.~\ref{fig:Unicell}),  we encode additional identifications of the triangles of the boundary: an edge of color 0 encodes the gluing of two triangles in the only way that respects the coloring (an edge of color $i\in\{1,2,3\}$ is glued to the other edge of color $i$, a vertex with two colors $i<j\in\{1,2,3\}$ is glued to the other vertex of colors $i,j$). Therefore, a 4-colored graph whose 0-residue (the 3-colored subgraph obtained by removing the edges of color 0) is dual to the boundary of the octahedron of Fig.~\ref{fig:Unicell} is itself dual to a self-gluing  (without boundary) of this octahedron. \\

More generally, if $G$ is a $(q+1)$-colored graph, each one of the connected components of the $q$-colored graph $G_{\hat 0}$ obtained by deleting all color-0 edges is dual to a $(q-1)$-dimensional colored triangulation, and the additional color-0 edges encode identifications of $(q-1)$-simplices in a unique way\footnote{In a $q$-dimensional triangulation, a vertex incident to a $q$-simplex can be associated the color of the only $(q-1)$-simplex it is not incident to. When gluing two $(q-1)$-simplices, the vertices that have the same colors are identified: this is done in a unique way. }. 
If we provide each $(q-1)$-dimensional triangulation with an interior, such spaces can therefore be considered as obtained from a collection of elementary $q$-dimensional building blocks with triangulated $(q-1)$-dimensional boundaries by identifying two-by-two the $(q-1)$-simplices of their boundaries.

The vertices of the colored graphs are dual to the $(q-1)$-simplices  that belong to the boundaries of the elementary building blocks. Moreover, the color-$0i$ cycles with $p$ color-0 edges are dual to the $(q-2)$-simplices with $p$ incident $(q-1)$-simplices in the dual discrete space. Therefore, the order $\delta_0$ of a $(q+1)$-edge-colored graph $G$ (given in 
Eq.~\eqref{eq:DeltaS0}) is a linear combination of the numbers $n_{q-2}$ and $n_{q-1}$ of $(q-2)$ and $(q-1)$-simplices of its dual space $\cM$,
\be 
\delta_0(G) = 1+\frac{q-1}2 n_{q-1}(\cM)-n_{q-2}(\cM).
\ee

Moreover, the 0-residues of a colored graph (the connected components of $G_{\hat 0}$) are dual to the elementary building blocks in the dual picture. 
Because of the additional connectivity condition, SYK graphs  encode gluings of a single building block, and are thus a generalization in higher dimensions  of unicellular maps  (discrete two-dimensional surfaces obtained from a single polygon by identifying two-by-two the edges of its boundary) \cite{unicell1, unicell2, unicell3}.

In dimension $q=2$, the 3-colored graphs are dual to bipartite combinatorial maps (which are orientable if and only if the colored graph is bipartite). A map is a drawing of a connected graph on a two-dimensional surface, considered up to homeomorphisms of the surface, and such that the connected components of the complement of the graph on the surface (called faces) are homeomorphic to discs (the map is said to be cellularly embedded). We only consider maps embedded on orientable surfaces in this paper. The order $\delta_0$ of a 3-colored graph $G$  is the excess of the dual bipartite map $\cM$ (or its number of independent cycles), 
\be 
\delta_0(G) = 1+E(\cM)-V(\cM).
\ee
Note that if the bipartite map $\cM$ is unicellular (i.e.~if it has a single face), then the order is twice the genus of the dual map, 
\be 
\delta_0(G) = 2g(\cM).
\ee
For unicellular spaces in dimension $q$, the order is therefore one possible generalization of the genus. 
One of the results we prove in this paper, is that for $q>2$, a large $(q+1)$-edge-colored graph of fixed order $\delta_0$ is almost surely dual to a unicellular colored space. This is not true in dimension 2: a large bipartite map of fixed excess is not a.a.s.~unicellular.


In dimension $q>2$, our results apply to the asymptotic enumeration of unicellular colored discrete spaces, according to a linear combination of their $(q-2)$ and $(q-1)$-dimensional elements, which reduces to twice the genus of the unicellular space for $q=2$. 

Although from a purely geometric point of view,  the choice of this particular combination of the number of $(q-1)$ and $(q-2)$ simplices appears to be arbitrary, the classification according to the order turns out to be remarkably simple with respect to other similar situations (for instance classifying colored triangulations or gluings of octahedra according to other linear combinations of the numbers of simplices, as done e.g.~in \cite{gilles, BLR, octa, LL}).

\medskip 

\subsection{Statement of the main results}\label{sec:state_main}

We let $\tilde g_{n,\delta}$ (resp.~$g_{n,\delta}$) be the number of rooted (resp.~rooted bipartite) colored graphs of fixed order $\delta$ with $2n$ vertices, and  $\tilde c_{n,\delta}$ (resp.~$c_{n,\delta}$) be the number of rooted (resp.~rooted bipartite) SYK graphs of fixed order $\delta$ with $2n$ vertices. In this paper, we will generally use analytic combinatorics tools and notations that can be found for example in \cite{flajolet}.

We also denote by $m_{\delta}$ the number of rooted\footnote{A map is called rooted if it has a marked edge that is given a direction.} trivalent maps with $2\delta-2$ vertices, which is given by  (see A062980 in OEIS) the recurrence
\[
m_1=1,\ m_2=5,\ m_{\delta}=(6\delta-8)m_{\delta-1}+\sum_{k=1}^{\delta-1} m_k\cdot m_{\delta-k}\ \ \ \mathrm{for}\ \delta\geq 3,
\]
and define the constant $\kappa_{\delta}$ as $\kappa_0=\sqrt{q/2\pi (q-1)^3}$
and 
\[
\kappa_{\delta}=\frac{1}{\Gamma\Big(\tfrac{3\delta-1}{2}\Big)}\frac 2 {q(q-1)}  \Bigl(\frac{q-1}{2q^3}\Bigr)^{{(3\delta-1)/}2}\Bigl(  \frac {q^4} 4 \Bigr)^{\delta}m_{\delta}, \ \ \mathrm{for}\ \delta\geq 1.
\]

\begin{theorem}\label{theo:counting} For $\delta\geq 0$ and $q\geq 3$, the numbers $g_{n,\delta}$,  $c_{n,\delta}$, $\frac 1 {2^\delta}\tilde g_{n,\delta}$, and $\frac 1 {2^\delta}\tilde c_{n,\delta}$ behave asymptotically as 
\be 
 \kappa_\delta\cdot n^{3(\delta-1)/2} \cdot \gamma^{n}, 
\ee 
where $\gamma = \frac {(q+1)^{q+1}}{q^q}$, and $\kappa_{\delta}$ is defined above.
\end{theorem}

\

As a consequence, for $q\geq 3$, if we let $G_{n,\delta}$ be a   random rooted $(q+1)$-edge-colored graph of order $\delta$ and with $2n$ vertices, then for $\delta$ fixed
and $n\to\infty$: 
\be 
\bP\bigl(G_{n,\delta}\text{ is SYK}\bigr) \rightarrow 1, \qquad\text{and}\qquad \bP\bigl(G_{n,\delta}\text{ is bipartite}\ |\ G_{n,\delta}\text{ is SYK}\bigr) \rightarrow 2^{-\delta}.
\ee 

\

Let us mention that the second estimate is reminiscent of what happens for the random quadrangulation $Q_{n,g}$ of genus $g$  with $n$ faces:  
for every $g\geq 0$  one has $\bP\bigl(Q_{n,g}\text{ is bipartite}\bigr)\rightarrow 2^{-2g}$ as $n\to\infty$
 (see~\cite[Theorem 2]{ChC} which gives a more general statement and related references).  

Proving Thm.~\ref{theo:counting} will be the main focus of the article. A first important tool is a bijection (introduced in~\cite{LL}) between colored graphs and so-called \emph{constellations} 
(certain partially embedded graphs defined later below Fig.~\ref{fig:Ex-Stacked}) such that  the order of a colored graph corresponds  to the excess of the associated constellation.
We recall this bijection in Section~\ref{sec:Bijections} and  also explain how it can be adapted to the non-orientable setting, giving the simple relation $\tilde g_{n,\delta}=2^\delta\ g_{n,\delta}$
(which is also reflected above by the fact that $g_{n,\delta}$ and $\tilde g_{n,\delta}/2^\delta$
have the same asymptotic estimate).  
Then, in Subsection~\ref{sec:kernel},  
we use the so-called method of kernel extraction to obtain an explicit expression
for the generating function $\cG_{\delta}(z)=\sum_{n\ge 0} c_{n,\delta} z^n$  
of rooted bipartite colored graphs of fixed order $\delta$ (we say that $z$ is the variable dual to $n$). 
Singularity analysis of the obtained expression then gives us the asymptotic estimate
of $g_{n,\delta}$. 
In Section~\ref{sec:SYK_conditions} we  then  give sufficient conditions (in terms of the 
kernel decomposition) for a colored graph of order $\delta$ to be an SYK graph,
and deduce from it that asymptotically, almost all colored graphs of order $\delta$
are SYK graphs. 
This implies that $c_{n,\delta}$ has the same asymptotic estimate as $g_{n,\delta}$;
similar arguments in the non-orientable case 
ensure 
that $\tilde c_{n,\delta}$ has the same asymptotic estimate as $\tilde g_{n,\delta}$.

\

Our results on the asymptotic combinatorial structure of graphs of fixed order allow us to obtain in Section~\ref{sec:Topo} the asymptotic topology of the boundary of the building block in the triangulation interpretation of Section~\ref{sub:Geom-Interp}. We indeed show that  if $G$ is a large random bipartite $(q+1)$-edge-colored graph of fixed order~$\delta$, then almost surely, it is dual to a unicellular space, and the boundary of the building block triangulates a piecewise-linear manifold whose topology is that of a connected sum of $\delta$ orientable $\cS^{q-2}$ bundles over $\cS^1$, $\cS^1\times\cS^{q-2}$, where $\cS^D$ denotes the $D$-dimensional sphere. 
In other words, if we denote by $\cT$ the triangulation dual to $G_{\hat 0}$ and $\lvert\cT\rvert$ the topological space triangulated by $\cT$, then
\be 
\lvert\cT\rvert \sim_{PL} \#_{\delta} (\cS^1\times\cS^{q-2}),
\ee 
where by  $\#_{\delta} X$ we mean the connected sum of $\delta$ copies of the orientable manifold $X$, and by $\sim_{PL}$ we mean ``is PL-isomorphic to''. This implies however that $G$ a.a.s.~does not represent a piecewise-linear manifold. Note that in the case $q=3$, the result implies that a 4-colored graph of fixed order $\delta$ is a.a.s.~dual to a 3-dimensional unicellular space whose building block has a 2-dimensional boundary of genus $\delta$.

\

\section{Bijection with constellations}
\label{sec:Bijections}

We recall here a bijection introduced in \cite{LL} from colored graphs of order $\delta$
to so-called \emph{constellations}\footnote{In \cite{LL}, constellations were called \emph{stacked maps}, as a central quantity in that work was the sum of faces over a certain set of submaps. Here we use the term \emph{constellations}, as it is a common name in the combinatorics literature for (the dual) of these objects.} of excess $\delta$. Thanks to this bijection, 
computing the generating function of colored graphs of fixed order amounts to computing
the generating function of constellations of fixed excess (which can be done using kernel extraction, as we will show in Section~\ref{sec:kernel}). 
We also explain in Section~\ref{sec:bij_non_bipartite} how 
the bijection can be adapted 
to
 the non-bipartite case. 

\medskip

\subsection{Bijection in the bipartite case}
\label{sub:Bij-Bip-case}

Given a bipartite colored graph $G\in\bG^q$, we first orient all the edges from black to white. We then contract all the color-0 edges as shown below in \eqref{fig:ContractColor0}, so that the pairs of black and white vertices they link collapse into $2q$-valent vertices which have one outgoing and one ingoing edge of each color $i\in\llbracket 1,q\rrbracket$. 
\be
\begin{array}{c} \includegraphics[scale=0.45]{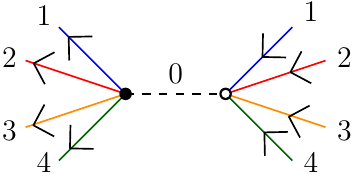} \end{array}
\qquad \longrightarrow \qquad \begin{array}{c} \includegraphics[scale=0.45]{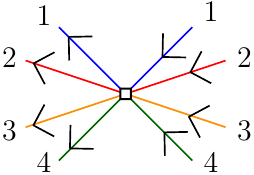} \end{array}
\label{fig:ContractColor0}
\ee 
The vertex resulting from the contraction of the distinguished color-0 edge is itself distinguished. The obtained (Eulerian) graph, $G_{/ 0}$, is such that the subgraph obtained by keeping only the color-$i$ edges ($i\in\llbracket 1,q\rrbracket$) is a collection of directed cycles. 
For each such color-$i$ cycle containing $p$ vertices, we add a color-$i$ vertex, and $p$ color-$i$ edges between that vertex and the $p$ vertices of the cycle,
and then we delete the original color-$i$ edges, as illustrated below.
\be
\label{fig:StarSubdivision}
\begin{array}{c} \includegraphics[scale=0.4]{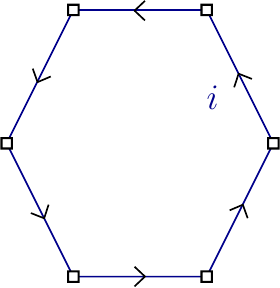} \end{array}
\qquad \longrightarrow \qquad \begin{array}{c} \includegraphics[scale=0.4]{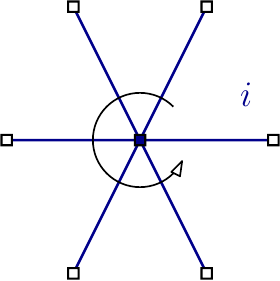} \end{array}
\ee 
 The cyclic ordering of the $p$ edges around the cycle translates into a  cyclic counterclockwise ordering of the $p$ edges around the color-$i$ vertex, 
and each of these (deleted) edges corresponds to a corner of the color-$i$ vertex. We say that a vertex is embedded if a cyclic ordering of its incident edges is specified.

\begin{figure}[h!]
\includegraphics[scale=0.7]{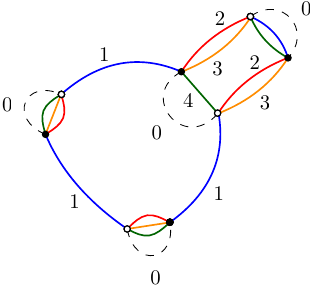}\hspace{2cm}\includegraphics[scale=0.7]{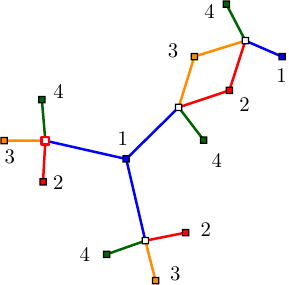}
\caption{\label{fig:Ex-Stacked}A bipartite SYK 5-colored graph and the corresponding  $4$-constellation.}
\end{figure}

Doing this operation at every color-$i$ cycle, we obtain a connected diagram $S=\Psi(G)$
 having
\begin{itemize}
\item non-embedded white vertices of valency $q$, with one incident edge of each color $i\in\llbracket 1,q\rrbracket$,
\item embedded color-$i$ vertices, 
\item color-$i$ edges, which connect a white vertex to a color-$i$ vertex,
\item one distinguished white vertex (resulting from the contraction of the root-edge).  
\end{itemize}

We denote by $\bS^q$ the set of such diagrams, which we call (rooted) $q$-constellations\footnote{We stress that in the usual definition of constellations, the white vertices are embedded, and the cyclic ordering of the edges is given by the ordering of their colors. 
However, we will not need here to view constellations as equipped with this canonical embedding.}.
An example is shown in Fig.~\ref{fig:Ex-Stacked}.
We recall that the excess of a connected graph $G$ is defined as $L(G) = E(G) - V(G) + 1$; it corresponds to its number of independent cycles. We let $\bS^q_{n,\delta}$ be the set of 
constellations in $\bS^q$ with $n$ white vertices and excess $\delta$, and let 
$\bG^q_{n,\delta}$ be the set of rooted bipartite $(q+1)$-edge-colored graphs with $2n$ vertices and order $\delta$. 

\begin{theorem} \cite{LL} \label{theo:bij}
The map $\Psi$ described above gives a bijection between $\bG^q_{n,\delta}$ 
and $\bS^q_{n,\delta}$, for every $q\geq 2$, $n\geq 1$ and $\delta\geq 0$.  
\end{theorem} 
 \proof 
The map  $\Psi$ is clearly invertible, hence gives a bijection from $\bG^q$ to $\bS^q$. 
Regarding the parameter
correspondence, for $S=\Psi(G)$, the color-$0$ edges of $G$ correspond to the white vertices
 of $S$, and these edges form a perfect matching in $G$, hence if $S$ has $n$ white
vertices then $G$ has $2n$ vertices. Finally, note that $\delta_0(G)=1+(q-1)n-F_0(G)$,
and $L(S)=E-n-m+1$, with $E$ the number of edges and $m$ the total number of colored
vertices in $S$.  Since each color-$0i$ cycle of $G$ is mapped to a color-$i$ vertex of $S$, we have $F_0(G)=m$. Since each edge of $S$ has exactly one white extremity and 
since white vertices have valency $q$, we have $E=qn$. Hence $L(S)=\delta_0(G)$. \qed

\ 

\begin{figure}[h!]
\raisebox{20.5pt}{\includegraphics[scale=0.9]{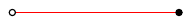}
 \hspace{1cm }$\longrightarrow $}  \hspace{1cm }\includegraphics[scale=0.9]{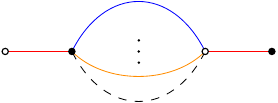} 
 \caption{Insertion of a pair of vertices linked by $q$ edges.}
\label{fig:InsertMelon}
\end{figure}

\begin{figure}[!h]
\raisebox{12pt}{\includegraphics[scale=1]{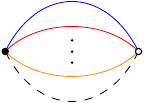}} \hspace{1.5cm }\includegraphics[scale=0.75]{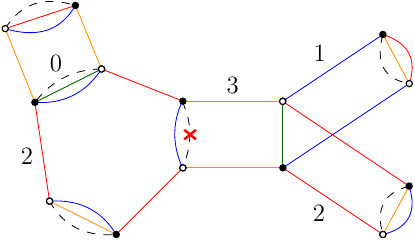}
\caption{\label{fig:MelonicGraphs} Melonic graphs.}
\end{figure}

\noindent{\bf Trees and melonic graphs.} Let us describe the colored graphs of vanishing order, which have been extensively studied in the literature. {\it Melonic graphs} are series-parallel colored graphs which appear in the context of random tensor models \cite{CritBehavior,Uncolor, Dartois}. 
They are obtained by recursively inserting pairs of vertices linked by $q$ edges,  as shown in Fig.~\ref{fig:InsertMelon}, starting from the only colored graph with two vertices (left of Fig.~\ref{fig:MelonicGraphs}).  Melonic graphs can equivalently be defined as the colored graphs in $\tilde \bG^q$ that satisfy the following identity:\footnote{They have a vanishing Gurau degree \cite{Uncolor}.}
\be 
q+(q-1)\bigl(\frac{1}2V(G) - R_0(G)\bigr) - F_0(G) = 0,
\ee 
or equivalently, 
\be 
\deltaS_0(G) = (q-1)\bigl(R_0(G) - 1\bigr).
\ee 
The colored graphs of order $\deltaS_0 = 0$,  i.e.~those that $\Psi$ maps to trees (constellations of vanishing excess), are melonic graphs which in addition have a single 0-residue (they are SYK graphs) \cite{blt}. Indeed, in the recursive construction of Fig.~\ref{fig:InsertMelon} for an SYK melonic graph, the edge on which a pair of vertices is inserted must not be of color 0, and it is easily seen that in a $q$-constellation, a white vertex incident to $q-1$ leaves whose colors are not $i\in\llbracket 1, q\rrbracket$ corresponds to an insertion as in Fig.~\ref{fig:InsertMelon} on an edge of color $i\neq 0$ in the colored graph (a leaf is a vertex of valency one).
\begin{prop}
\label{MelonsAreTrees}
The bijection $\Psi$ maps the melonic SYK graphs in $\bG^q$ to the trees in $\bS^q$.
\end{prop}
An example of a rooted melonic SYK graph is  shown on the right of Fig.~\ref{fig:MelonicGraphs} for $q=3$.

Note that this is also true for $q=2$, although in that case a melonic SYK 3-colored graph is dual to a unicellular planar bipartite map. Another remark is that the 0-residue $G_{\hat 0}$ is also a melonic graph, and in general, deleting the edges of a given color in a $(q+1)$-colored melonic graph, one is left with a collection of melonic graphs with $q$ colors.

\medskip

\subsection{The non-bipartite case}\label{sec:bij_non_bipartite}

Consider a (non-necessarily bipartite) $(q+1)$-edge-colored graph $G\in \tilde \bG^q$. It has an even number of vertices, since color-0 edges form a perfect matching. We assign an orientation to each non-root color-0 edge. If $G$ has $2n$ vertices, there are $2^{n-1}$ ways of doing so. 
A vertex is called an in-vertex (resp.~out-vertex) if it is the origin (resp.~end) of its incident color-0 edge. We then orient canonically the remaining half-edges (those on 
color-$i$ edges for $i\in \llbracket 1,q\rrbracket$); those at out-vertices are oriented outward 
and those at in-vertices are oriented inward. 

Contracting the color-0 edges into white vertices as in \eqref{fig:ContractColor0}, we obtain an Eulerian graph such that  for each $i\in \llbracket 1,q\rrbracket$,
 every vertex has exactly one ingoing half-edge and one outgoing half-edge of color $i$. 
We choose an arbitrary orientation for each cycle of color $i\in\llbracket 1,q \rrbracket$. For each white vertex and for each color $i$, the orientation of its incident half-edges either coincides with the orientation of the color-$i$ cycle it belongs to, or they are opposite. We perform a star subdivision as in \eqref{fig:StarSubdivision}, with the difference that now, each newly added color-$i$ edge carries a $\pm$ sign, $+$ if the orientations of the color-$i$ half-edges agree with that of the color-$i$ cycle, $-$ otherwise, as illustrated below.

\be
\label{fig:StarSubdivisionNonBip}
\begin{array}{c} \includegraphics[scale=0.4]{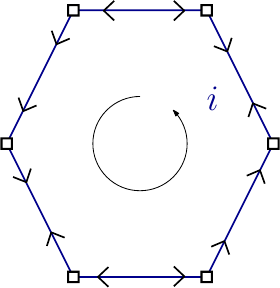} \end{array}
\qquad \longrightarrow \qquad \begin{array}{c} \includegraphics[scale=0.4]{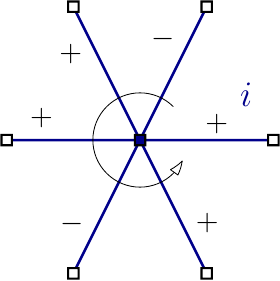} \end{array}
\ee 

We call \emph{signed colored graph} a rooted colored graph together with a choice of orientation of each non-root color-0 edge, and a choice of orientation of each color-$0i$ cycle for $i\in\llbracket 1,q\rrbracket$. We call \emph{signed constellation} a constellation (with a distinguished white vertex) 
 together with a choice of $\pm$ sign for every edge.
The transformation above defines a bijection $\hPsi$ between the set 
$\widehat\bG^q_{n, \delta}$ of signed colored graphs of order $\delta$ with $2n$ vertices
and the set $\widehat\bS^q_{n, \delta}$ of signed constellations with $n$ white vertices
and excess $\delta$. Let $\widetilde\bG^q_{n, \delta}$ be the set of rooted colored
graphs of order $\delta$ with $2n$ vertices. Since a colored graph $G\in \widetilde\bG^q_{n, \delta}$ 
 has $n-1$ non-root color-0 edges and satisfies $F_0(G)=1+(q-1)n-\delta$,
we have 
$\widehat\bG^q_{n, \delta}\simeq 2^{qn-\delta}\widetilde\bG^q_{n, \delta}$, where the symbol $\simeq$ means ``is isomorphic to".
Furthermore,
since a constellation in $\bS^q_{n, \delta}$ has $qn$ edges, we have 
$\widehat\bS^q_{n, \delta}\simeq 2^{qn}\times \bS^q_{n, \delta}$. Hence we obtain 
\be 
2^{qn-\delta} \times \widetilde\bG^q_{n, \delta}\  \simeq_{\hPsi}\  2^{qn} \times
 \bS^q_{n, \delta}\ \simeq_{\Psi}\ 2^{qn} \times \bG^q_{n, \delta}.
\ee

As a consequence, the  generating function $\widetilde\GF_{\delta}(z)$ of non-necessarily bipartite rooted $(q+1)$-edge-colored graphs of order $\delta$, with $z$ dual to half the number of vertices, satisfies:
\be 
\label{eq:Gen-Func-orient-vs-nonOrient}
\widetilde\GF_{\delta}(z) = 2^{\delta} \GF_{\delta}(z).
\ee 
We can 
thus focus on the bipartite case when dealing with the enumeration
of colored graphs of fixed order.

\medskip 

\section{Enumeration of colored graphs of fixed order}
In this section, we compute the generating function $\cG_\delta(z)$ 
of bipartite $(q+1)$-colored graphs of fixed order $\delta$. 
By Thm.~\ref{theo:bij} this is the 
generating function of constellations of excess $\delta$, with $z$ dual to the number of white
vertices. We rely on the so-called method of kernel extraction to obtain an explicit
expression of $\cG_\delta(z)$. Then singularity analysis on this expression will allow
us (in Subsection~\ref{sec:sing_analysis}) 
to obtain the asymptotic estimate of $g_{n,\delta}=[z^n]\cG_{\delta}(z)$,
 stated in Thm.~\ref{theo:counting}.  We recall that $g_{n,\delta}$ is the number of rooted colored graphs of fixed order $\delta$ with $2n$ vertices (we use the notation $[z^n]\cG_{\delta}(z)$ for the coefficient of the degree $n$ monomial of the generating function $G_\delta(z)$).
 
 \medskip

\subsection{Exact enumeration}\label{sec:kernel}

\begin{figure}[h!]
\includegraphics[scale=0.6]{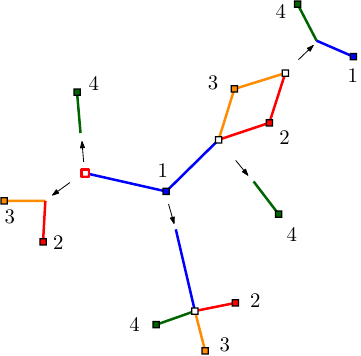}
\hspace{0.8cm}
\raisebox{30pt}{\includegraphics[scale=0.9]{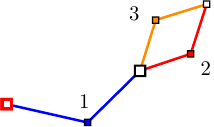}}
\hspace{1.2cm}
\raisebox{20pt}{\includegraphics[scale=1]{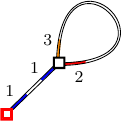}}
\caption{Cutting out tree contributions in the example of Fig.~\ref{fig:Ex-Stacked} (left) leads to its core diagram (center). The corresponding kernel is shown on the right of the figure.}
\label{fig:Pruning}
\end{figure}

For a constellation $S$, the \emph{core} $C$ of $S$ is obtained by iteratively deleting 
the non-root leaves (and incident edges) until all non-root vertices have valency at least $2$. 
This procedure is shown in Fig.~\ref{fig:Pruning} for the example of Fig.~\ref{fig:Ex-Stacked}. The core diagrams satisfy the following 
properties:
\begin{itemize}
\item
white vertices are non-embedded while $i$-colored vertices (for $i\in\llbracket 1,q \rrbracket$) are embedded, 
\item
white vertices have valency at most $q$, with incident edges of different colors,
\item
all non-root vertices (white or colored) have valency at least $2$,
\item
each edge carries a color $i\in\llbracket 1,q \rrbracket$, and connects a white vertex to a color-$i$ vertex.
\end{itemize}

We now focus on the maximal sequences of non-root valency-two vertices:

\begin{definition} A {\it chain-vertex} of a core diagram is a non-root vertex of valency two.
A {\it core-chain} is a path whose internal vertices are chain-vertices, but whose extremities are not chain-vertices. The {\it type} of a core-chain is given by the colors of its two extremities (colored or white), and by the color of their incident half-edges in the chain.
\end{definition}

Replacing all the core-chains by edges whose two half-edges retain the colors of the 
extremal edges on each side of the chain, we obtain the \emph{kernel} $K$ of the constellation $S$.
Note that $K$ is a diagram that has a distinguished white vertex (the root-vertex)
and satisfies the following conditions:
\begin{itemize}
\item
white vertices are non-embedded while $i$-colored vertices (for $i\in\llbracket 1,q \rrbracket$) are embedded, 
\item
white vertices have valency at most $q$, with incident half-edges of different colors,
\item
all non-root vertices (white or colored) have valency at least $3$,
\item
each half-edge carries a color $i\in\llbracket 1,q \rrbracket$, and is incident either to a white vertex
or to a vertex of color $i$. 
\end{itemize}
We call \emph{kernel diagrams} the (connected) graphs satisfying these properties. The excess of such 
a diagram $K$ is as usual defined as $E-V+1$, with $E$ its number of edges and $V$
its number of vertices. 
An important property is that a constellation and its kernel have  
equal excess. 
We let $\bK^q$ be the family of kernel diagrams and $\bK^q_{\delta}$ 
those of excess $\delta$. 
Since every non-root vertex in a kernel diagram has valency at 
least $3$, it is an easy exercise to show that $K$ has at most $3\delta-1$ edges (this calculation
will however be detailed in Section~\ref{sec:sing_analysis}), so that $\bK^q_{\delta}$ is a finite set.

An edge of $K$ is called unicolored (resp.~bicolored) if 
its two half-edges have the same color (resp. have different colors). 
For $K\in\bK^q$ we let $\cVw(K)$ and $\cVb(K)$ be the sets of white vertices
and of colored vertices in $K$, and let $\Vw=\mathrm{Card}(\cVw)$ and $\Vb=\mathrm{Card}(\cVb)$;  
we let $E(K),\Esame(K),\Edif(K)$ be the numbers of edges, of 
unicolored edges, and of bicolored edges in $K$; we also use refined notation  
$\Ebb(K)$, $\Ebbsame(K)$, $\Ebbdif(K)$, $\Ebw(K)$, $\Ebwsame(K)$, $\Ebwdif(K)$, $\Eww(K)$, $\Ewwsame(K)$, $\Ewwdif(K)$ for  the numbers of any/unicolored/bicolored edges whose extremities are colored/colored (resp.~colored/white, 
resp. white/white).  

\begin{figure}[h!]
\includegraphics[scale=0.6]{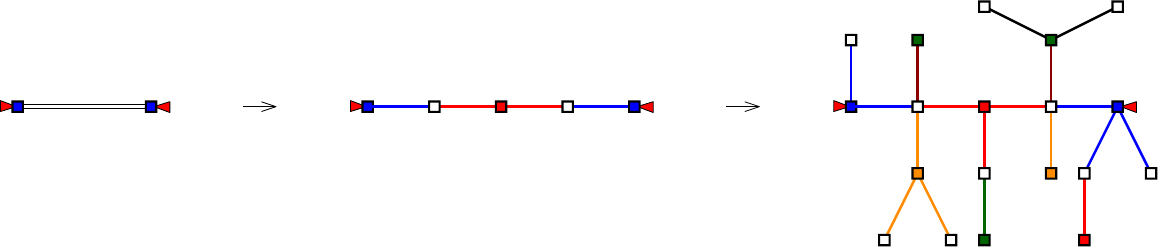}
\caption{\label{fig:RecoveringMaps} The constellations in $\bS_{\delta,K}^q$ are obtained from $K$ by replacing the edges by sequences of valency-two vertices and then inserting trees in the corners.}
\end{figure}

For $K\in\bK^q_{\delta}$ we let $\bS^q_{\delta,K}$
be the set of $q$-constellations (all of excess $\delta$) that have $K$ as kernel.  
A constellation in $\bS^q_{\delta,K}$ is generically obtained from $K$ where each edge $e$ is replaced by a core chain of the right type, i.e.~a sequence of valency-two vertices of arbitrary length, alternatively colored and white, which respects the boundary conditions, in the sense that extremal edges match the colors of the two half-edges that compose $e$, and an extremal vertex of the chain is white if and only if the incident extremity is white. Colored leaves are then added to white vertices so that they have  one incident edge of each color $i\in\llbracket 1,q\rrbracket$. An arbitrary tree rooted at a color-$i$ corner is then inserted at every color-$i$ corner (Fig.~\ref{fig:RecoveringMaps}). 

To obtain the generating function $\cG_{\delta,K}(z)$ of the family $\bS^q_{\delta,K}$, with $z$ dual to the number of white vertices, one must therefore take the product of the generating functions of the core-chains whose types correspond to the coloring of vertices and of half-edges in $K$, together with a certain number of tree generating functions. The generating function $\GF_T(z)$ of $q$-colored tree constellations rooted on a color-$i$ corner (for any fixed $i\in\llbracket 1,q\rrbracket$) and counted according to their number of white vertices is given by
 \be
 \GF_T(z)= 1+z  \GF_T(z)^q.
\ee
Its coefficients are the Fuss-Catalan numbers: $[z^n]\GF_T(z)=\frac 1{qk+1}\binom{qk+1}{k}$.

\begin{prop} 
\label{prop:ExEnum}
For $q\geq 3$, the generating function of rooted bipartite 
$(q+1)$-edge-colored graphs of fixed order $\delta$ is expressed
as $\cG_\delta(z)=\sum_{K\in\bK^q_{\delta}}\cG_{\delta,K}(z)$, where for each $K\in\bK^q_\delta$,
\be 
\label{eq:ExEnum}
\GF_{\delta,K}(z) = \biggl[\cB_\cK(y)\Bigl[\frac{y}{\bigl(1+y \bigr)\bigl(1-(q-1)y\bigr)}\Bigr]^{E(\cK)} y^{V_\circ(\cK)}\biggr]_{y=z\GF_T(z)^{q}},
\ee 
with the notation 
\be 
\label{eq:FactorExactEnum}
\cB_K(y)=[(q-1)y]^{\Ebb^=(K)} 
[1/y-q+2]^{\Ebw^=(K)+\Eww^=(K)}.
\ee 
\end{prop}

\proof
\
Following the approach  of
 \cite{LL}, we first compute the generating functions of 
core-chains of various types, in two variables $\zw,\zb$, where $\zw$ (resp.~$\zb$)  
is dual to the number of non-extremal white (resp.~colored) vertices in the chain.  
We let $y=z_\circ z_\bullet$. 
For $i,j\in\llbracket 1,q \rrbracket$
we let $\GFbb^{ij}(\zw,\zb)$ be the generating function of core-chains whose extremal vertices are colored and extremal edges have colors $i,j$ respectively. 
By symmetry of the role played by the colors, 
for every $i\neq j$ the generating functions $\GFbb^{ij}(\zw,\zb)$ 
are all equal to a common generating
function denoted by $\GFbb^{\neq}(\zw,\zb)$, and for every $i\in\llbracket 1,q \rrbracket$ the generating functions 
$\GFbb^{ii}(\zw,\zb)$ are all equal to a common generating function denoted by $\GFbb^{=}(\zw,\zb)$. 
A decomposition by removal of the first white vertex along the chain yields the system
\[
\GFbb^{\neq}(\zw,\zb)=\zw+(q-2)y\GFbb^{\neq}(\zb,\zw)+y\GFbb^=(\zb,\zw),\ \ \ \GFbb^=(\zb,\zw)=(q-1)y\GFbb^{\neq}(\zw,\zb),
\]
whose solution is
\be\label{eq:ker_1}
\GFbb^{\neq}(\zw,\zb)=\frac{\zw}{(1+y)(1-(q-1)y)},\ \ \GFbb^{=}(\zb,\zw)=(q-1)y\GFbb^{\neq}(\zb,\zw).
\ee
Similarly, we use the notations $\GFbw^{ij}(\zw,\zb), \GFbw^{=}(\zw,\zb), \GFbw^{\neq}(\zw,\zb)$
for the generating functions of core-chains whose extremal vertices are colored/white. 
By deletion of the extremal white vertex we find
\be \label{eq:ker_2}
\GFbw^{\neq}(\zw,\zb)=\zb\GFbb^{\neq}(\zw,\zb),\ \ \GFbw^{=}(\zw,\zb)=1+\zb\GFbb^{=}(\zw,\zb).
\ee
Finally, we use the notations $\GFww^{ij}(\zb,\zw), \GFww^{=}(\zb,\zw), \GFww^{\neq}(\zb,\zw)$
for the generating functions of core-chains whose extremal vertices are white/white.
By deletion of the extremal white vertices we find
\be \label{eq:ker_3}
\GFww^{\neq}(\zw,\zb)=\zb^2\GFbb^{\neq}(\zw,\zb),\ \ \GFww^{=}(\zw,\zb)=\zb(1+\zb\GFbb^{=}(\zw,\zb)).
\ee

Given a kernel diagram $\cK\in\bK^q$, let
\[
\cA_K(\zw,\zb):=(\GFbb^=)^{\Ebb^=(K)}(\GFbb^{\neq})^{\Ebb^{\neq}(K)}(\GFbw^=)^{\Ebw^=(K)}(\GFbw^{\neq})^{\Ebw^{\neq}(K)}(\GFww^=)^{\Eww^=(K)}(\GFww^{\neq})^{\Eww^{\neq}(K)}.
\]

Using the expressions \eqref{eq:ker_1}-\eqref{eq:ker_3}, 
we find

\begin{align*} 
\cA_\cK(z_\circ, z_\bullet)&= z_\bullet^{2E_{\circ\circ}(\cK) + E_{\circ\bullet}(\cK)} 
\cB_\cK( y )
 \Bigl[\frac{z_\circ}{\bigl(1+y \bigr)\bigl(1-(q-1)y\bigr)}\Bigr]^{E(\cK)}\\
&= z_\bullet^{\sum_{v_\circ \in \cVw(K)} \val(v_\circ)}  
\cB_\cK( y )
 \Bigl[\frac{z_\circ}{\bigl(1+y \bigr)\bigl(1-(q-1)y\bigr)}\Bigr]^{E(\cK)},
\end{align*}
where  $\val(v_\circ)$ is the valency of the vertex $v_\circ$, and where
\[ 
\cB_K( y )=[(q-1)y]^{\Ebb^{=}(K)} 
[1/y -q+2]^{\Ebw^{=}(\cK)+\Eww^{=}(K) }.
\]  
We then obtain
\[ 
\GF_{\delta,K}(z) = \cA_K\bigl(z \GF_T(z)^{q-2}, \GF_T(z)^{2}\bigr) \prod_{v_\circ\in \cVw(K)} z\GF_T(z)^{q-\Val(v_\circ)}\prod_{v_\bullet\in \cVb(K)} \GF_T(z)^{\Val(v_\bullet)},
\] 
where $z_\circ$ (resp.~$z_\bullet$) has been replaced by $z \GF_T(z)^{q-2}$ (resp.~$ \GF_T(z)^{2}$),
to account for the tree attachments. 
This rearranges into
\be 
\nonumber
\GF_{\delta,K}(z) = \biggl[ \cB_K(y) \Bigl[\frac{y}{\bigl(1+y \bigr)\bigl(1-(q-1)y\bigr)}\Bigr]^{E(K)} y^{V_\circ(\cK)}\biggr]_{y=z\GF_T(z)^{q}}.
\ee 
\qed

\subsection{Singularity analysis}\label{sec:sing_analysis}

We can now obtain 
the singular expansion of $\cG_{\delta,K}(z)$ for every given $K\in\bK^q_\delta$,
which yields the singular expansion of $\cG_{\delta}(z)$ and the asymptotic estimate of 
$g_{n,\delta}$ stated in Thm.~\ref{theo:counting}. 

We start with the singularity expansion of the tree generating function $\cG_T(z)$. 
From the equation $\cG_T(z)=1+z\cG_t(z)^q$ it is easy to find (see~\cite{CritBehavior}) that the dominant 
singularity of $\cG_T(z)$ is 
\be 
\label{eq:SingTrees}
z_c = \frac {(q-1)^{q-1}}{q^q}, \qquad \mathrm{with}\ \ \ \GF_T(z_c) = \frac q {q-1},
\ee
and we have the singular expansion 
\be  
\label{eq:SingGenTrees}
\GF_T(z) = \frac q {q-1} - \sqrt{\frac{2q}{(q-1)^3} \bigl( 1-\frac z {z_c}\bigr)}  + o\bigl( \sqrt{z_c-z}\bigr).
\ee

Using \eqref{eq:SingTrees}, we have $z_c\GF_T(z_c)^q=\GF_T(z_c)-1=1/(q-1)$, so that 
for any $K\in\bK_\delta$ we have
\be 
\GF_{\delta,K}(z) \sim \Bigl[\frac{1}{q\bigl(1-(q-1)z\GF_T(z)^{q}\bigr)}\Bigr]^{E(\cK)} \frac 1 {(q-1)^{V_\circ(\cK)}}.
\ee 
From \eqref{eq:SingTrees} and \eqref{eq:SingGenTrees}, we have
\be 
(q-1)z\GF_T(z)^q = 1 - \sqrt{\frac{2q}{q-1}\bigl( 1-\frac z {z_c}\bigr)}  + o\bigl( \sqrt{z_c-z}\bigr),
\ee 
and therefore, using the expression of Prop.~\ref{prop:ExEnum}, we find 
\be 
\GF_{\delta,K}(z) \sim \Bigl[\frac{q-1}{2q^3\bigl( 1-\frac z {z_c}\bigr)}\Bigr]^{\frac{E(\cK)}2} \frac 1 {(q-1)^{V_\circ(\cK)}}.
\ee 

The singularity exponent is thus maximal for kernel diagrams that have maximal number of edges (at fixed excess $\delta$). As a kernel diagram $K\in\bK_\delta^q$ has no vertices of valency one or two, apart maybe from the root, 
\[
2E(\cK) = \sum_{v\in \cK} \val(v) \ge 3(V(\cK) - 1) + 1 = 3V(\cK) - 2, 
\] 
with equality if and only if the root vertex has valency one, and all the other vertices have valency three. This implies that 
\[
\delta = E(\cK) - V(\cK) +1 \ge \frac 1 3 (E(\cK) + 1).
\] 
The maximal number of edges of a kernel diagram with fixed excess $\delta$ is therefore $3\delta - 1$, and we denote by $\boK^q_\delta$ the subset of those diagrams in $\bK^q_\delta$, i.e., kernel 
diagrams with a root-leaf and all the other vertices of valency $3$. We hence obtain
\be 
\GF_{\delta}(z) = \Bigl(\frac{q-1}{2q^3\bigl( 1-\frac z {z_c}\bigr)}\Bigr)^{\frac{3\delta-1}2}  \sum_{K\in\boK^q_\delta} 
\frac 1 {(q-1)^{V_\circ(\cK)}} + o\Bigl(\frac{1}{z_c-z}\Bigr)^{\frac{3\delta-1}2} .
\ee

We  orient cyclically the edges at each white vertex by the natural order of the colors they carry. For each non-root white vertex, there are $3 \binom q 3$ ways of choosing the colors of the incident half-edges, so that they are ordered correctly. Moreover, there are $q$ ways of choosing the color $i\in\llbracket 1,q \rrbracket$ of the half-edge incident to the root vertex, as well as the color of each colored trivalent vertex (this fixes the color of the incident half-edges).
Let $\bM_\delta$ be the set of maps with one leaf (called the root) 
and $2\delta-1$ other vertices all of valency $3$; note that these maps have $3\delta-1$ edges
hence have excess $\delta$. In addition, the cardinality of $\bM_\delta$ is equal to the coefficient $m_{\delta}$ introduced in Section~\ref{sec:state_main} (the number of rooted trivalent maps with $2\delta- 2$ vertices), as the trivalent vertex incident to the root leaf can canonically be replaced by an oriented edge.

From the preceding discussion we obtain
\be 
\sum_{K\in \boK^q_\delta}  \frac 1 {(q-1)^{V_\circ(\cK)}}
= \frac{q}{q-1}\sum_{M\in \bM_\delta}\left(q+\frac1{q-1}3\binom{q}{3} \right)^{2\delta-1}
= \frac{q}{q-1}m_{\delta}\Big(\frac{q^2}{2} \Big)^{2\delta-1},
\ee
where the factor raised to power $2\delta-1$ corresponds to the choice
 for each vertex of valency $3$ whether it is colored or white.

Finally, we obtain the following expression:
\begin{align*}
\GF_{\delta}(z) &= \frac q {(q-1)}  \Bigl(\frac{q-1}{2q^3\bigl( 1-\frac z {z_c}\bigr)}\Bigr)^{\frac{3\delta-1}2}\Bigl(  \frac {q^2} 2 \Bigr)^{2\delta - 1} 
m_\delta    + o\Bigl(\frac{1}{z_c-z}\Bigr)^{\frac{3\delta-1}2}\\
&= \frac 2 {q(q-1)}  \Bigl(\frac{q-1}{2q^3\bigl( 1-\frac z {z_c}\bigr)}\Bigr)^{\frac{3\delta-1}2}\Bigl(  \frac {q^4} 4 \Bigr)^{\delta} 
m_\delta    + o\Bigl(\frac{1}{z_c-z}\Bigr)^{\frac{3\delta-1}2}. 
\end{align*}

Using transfer theorems of singularity analysis~\cite[Corollary VI.1]{flajolet}, we conclude that, for $\delta\geq 1$,
\be
\label{eq:AsymptEstimate}
[z^n]\GF_{\delta}(z)=\frac{1}{\Gamma\Big(\tfrac{3\delta-1}{2}\Big)}\frac 2 {q(q-1)}  \Bigl(\frac{q-1}{2q^3}\Bigr)^{\frac{3\delta-1}2}\Bigl(  \frac {q^4} 4 \Bigr)^{\delta}m_{\delta}\cdot n^{3(\delta-1)/2}\cdot z_c^{-n}.
\ee
This gives the asymptotic estimate of the number $g_{n,\delta}=[z^n]\GF_{\delta}(z)$  of rooted bipartite colored graphs of fixed order $\delta$ with $2n$ vertices
in Thm.~\ref{theo:counting}. 

\medskip 

\section{The connectivity condition and SYK graphs}\label{sec:SYK_conditions}

In Subsections~\ref{subsec:ConnCond}, \ref{subsec:g-bigger-3} and \ref{subsec:g-is-3} we give sufficient conditions for a bipartite $(q+1)$-edge-colored graph to be an SYK
graph; we then deduce (using again singularity analysis) that the non-SYK graphs
have asymptotically negligible contributions, which ensures that in Thm.~\ref{theo:counting}, the coefficients 
$c_{n,\delta}$ have the same asymptotic estimate as the coefficients $g_{n,\delta}$ (the
latter estimate having been established in the last section). 
In Subsection~\ref{subsec:Conn-Non-Bip}, we adapt these arguments to non-necessarily bipartite graphs.

\subsection{Preliminary conditions}
\label{subsec:ConnCond}


\begin{definition} 
We say that a white vertex of a constellation $S\in\bS^q$ 
is {\it admissible}, if the corresponding two vertices in the colored graph $\Psi^{-1}(S)$ 
are linked in the graph by a path containing no edge of color~0. 
\end{definition}


\begin{claim}
\label{cor:Admiss}
A constellation is the image of an SYK graph if and only if all of its white vertices are admissible.
\end{claim}
\medskip

In both subsections below, we will need the following lemmas. We say that a vertex $v$ has a tree attached to it if one of the incident edges $e=(v, v')$ is a bridge (the number of connected components increases when the edge is removed), and the connected component that contains $v'$ after the removal is as tree.

\begin{lemma} 
\label{lemma:AdmVertexTree}
A white vertex with at least one tree attached 
to it
is admissible.
\end{lemma}
\proof Consider such a white vertex $v_\circ$ in a map $S\in\bS^q$, and a tree attached to it via  an edge of color $i\in\llbracket 1, q\rrbracket$. We prove the lemma recursively on the size of this tree. 

\begin{figure}[h!]
 \includegraphics[scale=1]{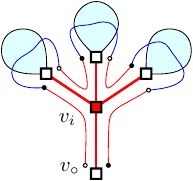} 
\hspace{2cm} \includegraphics[scale=1]{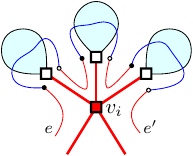} 
\caption{Concatenation of colored paths in Lemmas~\ref{lemma:AdmVertexTree} and \ref{lemma:ColoredPathTree}.}
\label{fig:ConcPAth}
\end{figure}

If the tree is just a color-$i$ leaf, it represents in the colored graph $G=\Psi^{-1}(S)$ an edge of color-$i$ between the corresponding two vertices in $G$, so that $v_\circ$ is admissible. If the tree has at least one white vertex, then the color-$i$ neighbor $v_i$ of $v_\circ$ has valency greater than one. All the other white vertices attached to $v_i$ have a smaller tree attached, and from the recursion hypothesis, they are admissible. To each corner of $v_i$ corresponds a color-$i$ edge in $G$, so that we can concatenate the colored paths in $G$ linking the pairs of vertices for all of these white vertices, as illustrated on the left of Fig.~\ref{fig:ConcPAth}. This  
concludes the proof.
\qed

\

Consider two color-$i$ edges $e$ and $e'$ in a colored graph $G$ corresponding to two corners incident to the same color-$i$ vertex $v_i$ in $S=\Psi(G)$, as shown on the right of Fig.~\ref{fig:ConcPAth}. These corners split the edges incident to $v_i$ into two sets $V_a$ and $V_b$.

\begin{lemma} 
\label{lemma:ColoredPathTree}
With these notations, if all the edges in either $V_a$ or $V_b$ lead to pending trees, then there exists a path in $G$ containing both $e$ and $e'$, without any color-$0$ edge. 
\end{lemma}

\proof 
Suppose that the condition of the lemma holds for $V_a$. All the white extremities of edges in $V_a$ have a tree attached, so that from Lemma~\ref{lemma:AdmVertexTree}, they are admissible. As above in the proof of Lemma~\ref{lemma:AdmVertexTree}, we can concatenate the colored paths for all these white vertices, using the color-$i$ edges incident to the corners between the edges in $V_a$, as shown on the right of Fig.~\ref{fig:ConcPAth}. 
\qed

\

Consider a $q$-constellation $S\in\bS^q$, its core diagram $C$, and its kernel diagram $K\in \bK^q$. We call tree contributions, the maximal trees removed to obtain the core diagram $C$ of $S$.  Consider a white vertex $v\in S$. We will say that it also belongs to $C$ if it is not internal to a tree contribution, and we will say that it also belongs to $K$ if,  in addition,  it is not a chain-vertex of $C$.

\begin{lemma} 
\label{lemma:VertValAdm}
With these notations, if $v$ is of valency smaller than $q$ in $K$, then it is admissible in $S$.
\end{lemma}
\proof If $v$ is of valency $d<q$ in $\cK$, it means that $q-d>0$ tree contributions have been removed in the procedure leading from a constellation $S$ to its kernel diagram $\cK$. We conclude applying Lemma~\ref{lemma:AdmVertexTree}. \qed

\medskip

\begin{lemma} 
\label{lemma:ValqResidues}
Let $G\in\bG^q$ and let $K$ be the kernel of the constellation associated to $G$. 
If $K$ has no white vertex of valency $q$, then $G$ is an SYK graph. 
\end{lemma}
\proof  Let $S\in \bS_{\delta,K}$. From Lemma~\ref{lemma:VertValAdm}, the vertices of $S$ that also belong to $K$ are admissible, as they are of valency smaller than $q$. 
The other white vertices of $S$ necessarily have a tree attached: either they are internal to a tree contribution, either they are chain-vertices of the corresponding core diagram. We conclude applying Lemma~\ref{lemma:AdmVertexTree} to every white vertex, and then the Claim~\ref{cor:Admiss}. \qed

\

\subsection{The case $q>3$}
\label{subsec:g-bigger-3}

For $\delta\geq 1$ and $q\geq 3$, a $(q+1)$-edge-colored graph is called \emph{dominant} 
if the kernel of its associated constellation 
belongs to 
 $\boK_\delta^q$, the set of kernel 
diagrams with a root-leaf and all the other vertices of valency $3$. 

\begin{coroll}
\label{cor:AsymptSYK_1}
For $q>3$ and $\delta\geq 1$, 
a dominant $(q+1)$-edge-colored graph is an SYK graph. Hence
$g_{n,\delta}\sim c_{n,\delta}$ as $n\to\infty$.  
\end{coroll}
\proof The associated kernel $K$ has a root-vertex of valency $1$ and 
and all its other vertices are of valency $3$. Hence all the vertices of $K$ have valency smaller than $q$. 
 Hence, from Lemma~\ref{lemma:ValqResidues}, $G$ is an SYK graph. 

We have seen in Section~\ref{sec:sing_analysis} that the non-dominant $(q+1)$-edge-colored graphs have asymptotically
a negligible contribution. Hence $g_{n,\delta}\sim c_{n,\delta}$.
\qed

\

\subsection{The case $q=3$}
\label{subsec:g-is-3}

It now remains to show that  
for $q=3$ we have $c_{n,\delta}\sim g_{n,\delta}$.   
Note that we cannot just apply Lemma~\ref{lemma:ValqResidues} as in the case $q>3$, since all non-root 
vertices of a kernel diagram $K\in\bK^3$ have valency at least $3=q$.

\begin{lemma}  
\label{lemma:AdmNonEmptyChain}
Let $S\in\bS^q$ be a $q$-constellation, with $C$ its core diagram and $K$ its kernel diagram. 
Let $v_\circ$ be a white vertex that belongs to $S$, $C$, and $K$.

If there is at least one core-chain in $C$ incident to $v_\circ$ and 
 containing at least one internal white vertex, then $v_\circ\in S$ is admissible.
\end{lemma}
\begin{figure}[h!]
 \includegraphics[scale=1]{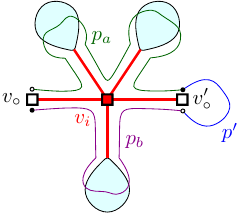}
 \caption{Concatenation of paths in the proof of Lemma~\ref{lemma:AdmNonEmptyChain}.}
 \label{fig:AdmNonEmptyChain}
\end{figure}

\proof 
Consider a core-chain in $C$ incident to $v_\circ$ and let $v'_\circ$ be the closest white chain-vertex in the chain (see Fig.~\ref{fig:AdmNonEmptyChain}). There is necessarily a color-$i$ chain-vertex for some $i\in\llbracket 1,q \rrbracket$ between $v_\circ$ and $v'_\circ$, which we denote by $v_i$ ($v_i$ is in the chain, at distance one from both $v_\circ$ and $v'_\circ$).  

The vertex  $v'_\circ$ has $q-2>0$ trees attached in $S$, and using Lemma~\ref{lemma:AdmVertexTree}, it is therefore admissible. We denote by $p'$ the corresponding path in $G=\Psi^{-1}(S)$. 

The vertex $v_i$ has two incident corners in $C$, both of which might have some tree contributions attached in $S$. These tree contributions are naturally organized in two groups $V_a$ and $V_b$ (which correspond to the two corners of $v'_\circ$ in $C$). Applying Lemma~\ref{lemma:ColoredPathTree} to both groups, we obtain two paths $p_a$ and $p_b$ in $G$. The concatenation of  $p'$, $p_a$ and $p_b$, gives a colored path between the two vertices corresponding to $v_\circ$ in $G$, so that $v_\circ$ is admissible.
\qed

\begin{lemma}  
\label{lemma:CSofSYK2}
Consider a  $(q+1)$-edge-colored graph $G$, and the core diagram $C$ of the $q$-constellation $S=\Psi(G)$. If every white vertex of $C$ either is of valency $d<q$, or has an incident core-chain containing at least one white chain-vertex, then $G$ is an SYK graph.
\end{lemma}

\proof This is a simple consequence of Lemma~\ref{lemma:VertValAdm}, Lemma~\ref{lemma:AdmNonEmptyChain}, and Claim~\ref{cor:Admiss}. \qed

\begin{lemma} 
\label{lemma:no_empty_white}
For $q\geq 3$ and $\delta\geq 1$, let $G$ be a  random edge-colored graph in $\bG^q_{n,\delta}$, 
and let $C$ be the core of the associated constellation $S$. Then, a.a.s.~all the 
core-chains of $S$ contain at least one (internal) white vertex. 
\end{lemma}
\proof
Let $r_{n,\delta}$ be the number of edge-colored graphs from $\bG^q_{n,\delta}$
with $n$ vertices, such that one of the core-chains is distinguished (i.e.~the kernel has a distinguished edge) with the condition 
that this distinguished core-chain has no internal white vertex. 
Lemma~\ref{lemma:CSofSYK2} ensures that $g_{n,\delta}-c_{n,\delta}\leq r_{n,\delta}$
hence we just have to show that $r_{n,\delta}=o(g_{n,\delta})$.
We let $\cR_\delta(z)=\sum_{n\geq 1} r_{n,\delta}z^n$ be the associated generating function. 
For every $K\in\bK_{\delta}$, the contribution to  $\cR_\delta(z)$ in the case where
the distinguished edge of $K$ has two white extremities and two half-edges of the same color
is (with the notations
in the proof of Prop.~\ref{prop:ExEnum}) equal to 
\[ 
\left[\frac{\Eww^=(K)\cdot(1+\zb)}{\GFww(\zw,\zb)}\cA_K\bigl(\zw,\zb)\right] \prod_{v_\circ\in \cVw(K)} z\GF_T(z)^{q-\Val(v_\circ)}\prod_{v_\bullet\in \cVb(K)} \GF_T(z)^{\Val(v_\bullet)},
\] 
where $\zw=\GF_T(z)^{q-2}$ and $\zb=\GF_T(z)^2$. It is then easy to check that, due
to the  $\GFww(\zw,\zb)$ appearing in the denominator, 
the leading term in the singular expansion is $O((z-z_c)^{-(E(K)-1)/2})$. 
This also holds for all the other possible types of the distinguished kernel edge, 
so that we conclude that $r_{n,\delta}=O(z_c^{-n}n^{(3\delta-4)/2})=o(g_{n,\delta})$.
\qed

\begin{theorem} 
\label{thm:AsymptSYK_2}
For $q\ge 3$ and $\delta\geq 1$, we have $c_{n,\delta}\sim g_{n,\delta}$ as $n\to\infty$. 
\end{theorem}
\proof
From Lemma~\ref{lemma:CSofSYK2} and Lemma~\ref{lemma:no_empty_white} 
it directly follows that (for $q\geq 3$ and $\delta\geq 1$) 
 the random  
edge-colored graph $G\in\bG^q_{\delta,n}$ is a.a.s.~an SYK graph. Hence
for $q\geq 3$ we have $c_{n,\delta}\sim g_{n,\delta}$. 
\qed

\medskip

\subsection{ The non-bipartite case}
\label{subsec:Conn-Non-Bip}

Let us go through the arguments of the last section, to adapt them in the case of generic colored graphs. 

Firstly,  choosing an orientation for every color-0 edge and color-$0i$ cycle does not change the number of 0-residues, so that we can work with signed colored graphs and signed constellations. Prop.~\ref{prop:Admiss}, Claim~\ref{cor:Admiss} are obviously true for signed colored graphs and signed constellations. Lemma~\ref{lemma:AdmVertexTree}, Lemma~\ref{lemma:VertValAdm} and Lemma~\ref{lemma:ValqResidues} also hold for signed constellations and signed colored graphs, as tree contributions represent bipartite (melonic) subgraphs. Therefore, Corollary~\ref{cor:AsymptSYK_1} is also valid for non-necessarily bipartite $(q+1)$-colored graphs, with $q>3$ and $\delta>0$.

Similarly, chains (and their attached trees) represent bipartite subgraphs of the colored graphs, so that for $q=3$, Lemmas~\ref{lemma:ColoredPathTree} and \ref{lemma:AdmNonEmptyChain} and \ref{lemma:CSofSYK2} can still be used without modification. It remains to adapt Lemma~\ref{lemma:no_empty_white} for signed colored graphs, i.e.~to prove that a.a.s, all core-chains of a signed constellation with $n$ vertices and excess $\delta$ contain at least one internal white vertex. This is true, as choosing a sign $\pm$ for every one of the $\delta + n -1$ edges does not modify this property.

Thus, Thm.~\ref{thm:AsymptSYK_2} generalizes for non-necessarily bipartite graphs in $\tilde G^q$: 
For $q\ge 3$ and $\delta\geq 1$, we have $\tilde c_{n,\delta}\sim \tilde g_{n,\delta}$ as $n\to\infty$. Thm.~\ref{theo:counting} follows from this, as well as \eqref{eq:Gen-Func-orient-vs-nonOrient} and \eqref{eq:AsymptEstimate}.

\medskip

\section{Topological structure of large random SYK graphs of fixed order}
\label{sec:Topo}

In this section, we consider a large colored graph $G\in \bG^q$ of fixed order $\deltaS_0(G) = \delta$ and that has a single residue, which we denote by $G_{\hat 0}$. In this case, the 0-residue can be interpreted as the graph dual to the triangulated boundary of the unique building block of the unicellular space dual to $G$ (see the introductory section \ref{sub:Geom-Interp} on colored triangulations). We further denote by $S=\Psi(G)$, $C$ its core diagram, and $\cK$ its kernel diagram, 
see again Fig.~\ref{fig:Pruning}.

\

We determine the a.a.s.~topology of the boundary of the building block of the large orientable unicellular space dual to the SYK graph $G$ of order $\delta$. This is done by adapting some results on the topology of triangulations dual to colored graphs \cite{Italians, GagliardiDipol, Handles,1N-arb, TensorTopology}, to the 0-residue $G_{\hat 0}$. In dimension $q=2$, an orientable unicellular map with $n$ edges is obtained from a disc whose polygonal boundary contains $2n$ edges. Thus, whatever the topology of the unicellular map (its genus), the topology of the boundary of the building block is always that of the circle $S^1$ (since the building block is always a disc). The topology of the building block is therefore a much weaker information than the topology of the glued space itself. Still, however, it provides some information about the possible spaces one can obtain from the building block. 

In Section~\ref{sub:Manifold}, we show that $G_{\hat 0}$ a.a.s.~represents a piecewise-linear manifold, and in Section~\ref{sub:AStopo}, we determine its a.a.s.~topology. In particular, we will see that although the boundary of the building block a.a.s.~triangulates a manifold, the glued space itself a.a.s.~never triangulates a piecewise-linear manifold (it contains a singularity and is therefore a pseudo-manifold). Still, fixing the order of the graph is responsible for the a.a.s.~non-singular topology of the building block (in \cite{carrance}, it is shown that a uniform $(q+1)$-colored graph with $q>2$ and all of its residues a.a.s.~represent pseudo-manifolds with singularities).

\medskip

\subsection{On the residues of $G_{\hat 0}$}
In order to obtain the results mentioned above, we need the following preliminary results.

We have seen in Lemma~\ref{lemma:no_empty_white} that, 
for $S$ a large random constellation of fixed order $\delta$,
then a.a.s.~every core-chain of $S$ has at least one white vertex. In view
of establishing the (a.a.s.)~topological structure of $S$, we need the following
complement to Lemma~\ref{lemma:no_empty_white} for colored vertices\footnote{With some more work
 it should be possible to show a limit law for the joint distribution of the lengths of the core-chains  
rescaled by $\sqrt{n}$, using generating functions and the method of moments.}.

\begin{lemma} 
\label{lemma:LongChainsColor}
Consider a large $(q+1)$-edge-colored graph $G$ of fixed order $\delta$, and the core diagram $C$ of the constellation $\Psi(G)$. Then a.a.s., any core-chain in $C$ contains at least one chain-vertex of every color  $i\in\llbracket 1, q \rrbracket$. 
\end{lemma}

\proof 
We proceed similarly as in the proof of Lemma~\ref{lemma:no_empty_white}.  
Let $a_{n,\delta}$ be the number of edge-colored graphs from $\bG^q_{n,\delta}$ where
all core-chains have at least one vertex in each color $i\in\llbracket 1, q \rrbracket$. And for $i\in\llbracket 1, q \rrbracket$ 
we let $b_{n,\delta,i}$ be the number of edge-colored graphs from $\bG^q_{n,\delta}$ with 
a distinguished kernel edge such  that  the corresponding core-chain has no vertex of color $i$. 
Clearly we have $g_{n,\delta}-a_{n,\delta}\leq \sum_{i=1}^q b_{n,\delta,i}$, so we 
just have to show that $b_{n,\delta,i}=o(g_{n,\delta})$ as $n\to\infty$. For $i\in\llbracket 1, q \rrbracket$ 
let $\cB_{\delta,i}(z)=\sum_{n\geq 1} b_{n,\delta,i}z^n$. 
For every $K\in\bK^q_\delta$ and $j\neq i$ the contribution to $\cB_{\delta,i}(z)$ where the kernel is $K$ 
and the distinguished edge of $K$ is white/white with both half-edges of color $j$ is equal to  
\[
\left[\Eww^{jj}(K)\frac{\cF_{\circ\circ}^{=}(\zw,\zb)}{\GFww^{=}(\zw,\zb)}\cA_K\bigl(\zw,\zb)\right] \prod_{v_\circ\in \cVw(K)} z\GF_T(z)^{q-\Val(v_\circ)}\prod_{v_\bullet\in \cVb(K)} \GF_T(z)^{\Val(v_\bullet)},
\]
where $\zw=\GF_T(z)^{q-2}$, $\zb=\GF_T(z)^2)$, $\cF_{\circ\circ}^{=}(\zw,\zb)$
 is the analog of $\GFww^{=}(\zw,\zb)$ with $(q-1)$ allowed colors  instead of $q$ colors 
 and $\Eww^{jj}(K)$ is the number of 
white/white edges of $K$ with both half-edges of color $j$. As we have already seen in Section~\ref{sec:sing_analysis}, we have
$\GFww^{=}(\GF_T(z)^{q-2},\GF_T(z)^2)=O((z-z_c)^{-1/2}$ as $z\to z_c$. On the other hand one can readily check that 
 $\cF_{\circ\circ}^{=}(\GF_T(z)^{q-2},\GF_T(z)^2)$ converges as $z\to z_c$  to a positive constant;
indeed the denominator of $\cF_{\circ\circ}^{=}(\GF_T(z)^{q-2},\GF_T(z)^2)$ involves
the quantity $(1-(q-2)z\cG_T(z)^q)$, which converges to $1/(q-1)$ as $z\to z_c$. 
Hence the above contribution  to $\cB_{\delta,i}(z)$  is $O((z-z_c)^{-(E(K)-1)/2})$. 
This also holds for all the other  contributions to $\cB_{\delta,i}(z)$,  
so that $b_{n,\delta}=o(g_{n,\delta})$. 
\qed

\

In the following, 
we specify the index $q$ for the bijection $\Psi_q : \bG^q \rightarrow \bS^q$ of Thm.~\ref{theo:bij}.
In order to characterize the topology of $G_{\hat 0}$, we need to study its $\C$-residues for $\C\in\llbracket 1, q\rrbracket$.  They are the connected components of the graph obtained from $G_{\hat 0}$ by deleting all the color-$\C$ edges. We denote by $S_{\hat \C}$ the (non-necessarily connected) constellation obtained from $S=\Psi_q(G)$ by deleting all the edges and vertices of color $\C$. The $\C$-residues of $G_{\hat 0}$ are subgraphs of the  $\C$-residues of $G$, which in turn are obtained by applying the inverse bijection $\Psi_{q-1}^{-1}$ to the connected components of $S_{\hat \C}$. We have the following simple consequence of Lemma~\ref{lemma:LongChainsColor}.

\begin{coroll} [of Lemma~\ref{lemma:LongChainsColor}]
\label{cor:iResidue}
With the notations of the lemma, for any $\C\in\llbracket 1, q\rrbracket$, the (non-necessarily connected) constellation $S_{\hat \C}$ is a.a.s.~a collection of trees.
\end{coroll}

\medskip

\subsection{$G_{\hat 0}$ asymptotically almost surely represents a manifold}
\label{sub:Manifold}

We recall that if $G$ is a $(q+1)$-colored graph, $\cT(G)$ its dual triangulation, and $\lvert \cT(G)\rvert$ the topological space that $\cT(G)$ triangulates, then $G$ is said to represent $\lvert \cT(G)\rvert$. We will use the following two topological results:

\begin{prop}[see e.g.~\cite{TensorTopology}, Prop.~3]
\label{prop:PLMan}
The $q$-colored graph $G_{\hat 0}$ represents a $(q-1)$-dimensional piecewise-linear (PL) manifold if and only if, for every color $\C\in\llbracket 1, q \rrbracket$, the $\C$-residues of  $G_{\hat 0}$  all represent $(q-2)$-spheres. 
\end{prop}

\begin{prop}[\cite{1N-arb}] 
\label{prop:TopoMelo}
Every $(q-1)$-colored melonic graph  
represents the $(q-2)$-sphere.
\end{prop}

We are precisely in the situation to use these two results:

\begin{lemma}
\label{lem:MelonResidu}
Consider a $(q+1)$-edge-colored graph $G$ with $S=\Psi_q(G)$, such that $S_{\hat \C}$ is a collection of trees. Then the $\C$-residues of $G_{\hat 0}$ are all melonic $(q-1)$-colored graphs.
\end{lemma}
\proof Suppose that $G$  is such that  the connected components of $S_{\hat \C}$ are trees.
The connected components of $G_{\hat \C}$ are obtained from the connected components of $S_{\hat \C}$ by applying the inverse bijection $\Psi^{-1}_{q-1}$. Using Prop.~\ref{MelonsAreTrees}, 
we know that $G_{\hat \C}$ is a collection of $q$-colored melonic graphs. This property is still satisfied when deleting all color-0 edges, thus every $\C$-residue of $G_{\hat{0}}$ is a melonic $(q-1)$-colored graph. 
\qed

\begin{lemma}
\label{lem:Topo1}
Consider a $(q+1)$-colored SYK graph $G$ with $S=\Psi_q(G)$, such that  for any $\C\in\llbracket 1, q\rrbracket$, $S_{\hat \C}$ is a collection of trees. Then $G_{\hat 0}$ represents a PL-manifold.
\end{lemma}
\proof This follows from Prop.~\ref{prop:PLMan}, Prop.~\ref{prop:TopoMelo} and Lemma~\ref{lem:MelonResidu}. \qed

\

We recall (Section~\ref{sub:Geom-Interp}), that by duality, an SYK graph is dual to a unicellular discrete space, obtained by gluing two-by-two the $(q-1)$-simplices  of its boundary. The boundary of its only building block is a colored triangulation dual to the connected $q$-colored graph $G_{\hat 0}$. 

\begin{theorem} 
\label{thm:Residue-Manifold}
Consider a large $(q+1)$-edge-colored graph $G$ of fixed order~$\delta$. Then almost surely, it is dual to a unicellular space, and the  boundary of the building block triangulates a PL-manifold.
\end{theorem}
\proof 
We have seen in Section~\ref{sec:SYK_conditions} that $G$ is a.a.s.~an SYK graph, which from Lemma~\ref{lemma:LongChainsColor} (Corollary~\ref{cor:iResidue}), a.a.s.~satisfies that for any $\C\in\llbracket 1, q\rrbracket$, $S_{\hat \C}$ is a collection of trees. Therefore, a large random $(q+1)$-edge-colored graph $G$ of fixed order~$\delta$ a.a.s.~satisfies the conditions (and thus the conclusions) of Lemma~\ref{lem:Topo1}. \qed

\

Note that this result applies both to bipartite and non-bipartite graphs. In addition to the conclusions above,  $G_{\hat0}$ represents an orientable manifold if and only if it is bipartite (see e.g.~\cite{TensorTopology}, p.~5).

\

\subsection{The a.a.s.~topology of $G_{\hat 0}$}
\label{sub:AStopo}

For large bipartite graphs $G\in \bG^q$ of fixed order, we are able to characterize in more details the a.a.s.~topology of the triangulation $\cT(G_{\hat 0})$ dual to $G_{\hat 0}$. We will need the following Lemma.

\begin{lemma}
\label{lem:Cutting-Trees}
Cutting off trees in a $q$-constellation $\Psi_q(G)$ does not change the topology of $\cT(G_{\hat 0})$.
\end{lemma}
\proof A non-empty tree contribution necessarily contains a white vertex attached to $q-1$ leaves. In the colored graph picture, this corresponds to a pair of vertices linked by $q$ colors, including 0. Two other external edges of a common color are attached to the two vertices. Deleting the pair of vertices in the colored graph and reconnecting the external pending half-edges amounts to deleting the white vertex and incident leaves in the constellation. In $G_{\hat 0}$, this is the same operation as in $G$ but with a pair of vertices that share $q-1$ edges. It is called a $(q-1)$-dipole removal and it is proven in \cite{GagliardiDipol} (Corollary 5.4) that the  $q$-colored graphs before and after a $(q-1)$-dipole removal represent the same topological space.\qed

\

If $S$ is a $q$-constellation, we denote by $S_{ij}$ the sub-diagram obtained from $S$ by keeping only the edges of color $i$ and $j$.

\begin{lemma} 
\label{lem:Forest-Bicolor}
With these notations, if $S_{ij}$ is a forest, the two vertices at the extremities of any color-0 edge in the colored graph $G=\Psi^{-1}(S)$,  belong to a common color-$ij$ cycle.
\end{lemma}
\proof Using Prop.~\ref{MelonsAreTrees}, 
by applying $\Psi_{2}^{-1}$, each connected component $s$  of $S_{ij}$ is mapped to an SYK  3-colored melonic graph $\Psi_{2}^{-1}(s)$, which is therefore connected when the color-0 edges are deleted: $\bigl[\Psi_{2}^{-1}(s)\bigr]_{\hat 0}$ is a connected color-$ij$ cycle  of $G_{\hat 0}$. The two vertices corresponding to a given white vertex of a connected component of $S_{ij}$ belong to the same connected 3-colored graph, and therefore to the same color-$ij$ cycle. \qed

\ 

Consider a constellation $S$, and its core diagram $C$.
A white chain-vertex in a core-chain has only two incident edges that belong to the core-chain (i.e.~that are not bridges). These edges have two different colors $i$ and $j$ in $\llbracket 1, q\rrbracket$. We denote by $v_{ij}$ such a chain-vertex.

\begin{definition} We say that a chain-vertex $v_{ij}$ in  $C\in \bC^{q}$ is a  {\it handle-vertex}  if the corresponding two vertices in $\Psi_q^{-1}(S)$ are linked in the graph by a color-$ij$ cycle, where $i,j\neq 0$.
\end{definition}

\begin{figure}[h!]
 \includegraphics[scale=1]{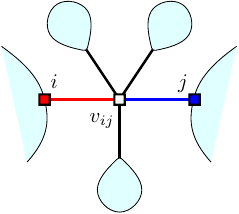} \hspace{1cm}\raisebox{12ex}{$\rightarrow$}\hspace{1cm}\raisebox{5.4ex}{\includegraphics[scale=1]{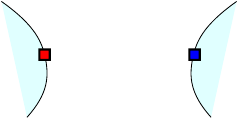} }
 \caption{Deletion of a chain-vertex and incident trees in a constellation.}
 \label{fig:Chain-vert-del}
\end{figure}

At this point, we must define an operation on the constellation $S$:  the deletion of a white chain-vertex $v_{ij}$, of the trees attached, and of the two incident edges of color $i$ and $j$. We denote by $S_{\setminus v_{ij}}$ the resulting constellation. This operation is illustrated in Fig.~\ref{fig:Chain-vert-del}.

We respectively denote by $\cT$ and $\cT_{\setminus v_{ij}}$ the two $(q-1)$-dimensional triangulations dual to   $(\Psi_q^{-1}(S))_{\hat 0}$ and   $(\Psi_q^{-1}(S_{\setminus v_{ij}}))_{\hat 0}$, and $\lvert\cT\rvert$ and $\lvert\cT_{\setminus v_{ij}}\rvert$ the topological spaces they triangulate. 
The following theorem is a translation  in the context of constellations of the main result of \cite{Handles}, Thm.~7, in the case of orientable spaces (see also the first point of Remark 2 in that reference). 

\begin{theorem}[\cite{Handles}]
\label{thm:Handles-Topo}
 If $v_{ij}$ is a handle-vertex, if $\lvert \cT\rvert$  is a connected PL-manifold, and if $\cT_{\setminus v_{ij}}$ is also connected, then $\lvert \cT\rvert$ is the connected sum of $\lvert \cT_{\setminus v_{ij}}\rvert$ and $\cS^1\times\cS^{q-2}$
 \be
\lvert \cT\rvert \sim_{PL} \lvert\cT_{\setminus v_{ij}}\rvert\,\#\, (\cS^1\times\cS^{q-2}).
 \ee
\end{theorem}

Note that in \cite{Handles}, the handles are defined in the $q$-colored graphs as pairs of vertices linked by $q-2$ colors, leaving two incident edges of some color $i$ and two of some color $j$, which all belong to a  common color-$ij$ cycle.
In the constellation picture, this means that the handle-vertex is incident to $q-2$ leaves. However, from Lemma~\ref{lem:Cutting-Trees}, removing or adding the tree contributions does not change the topology, so we chose to state the theorem in this more general setting.

\begin{lemma}
\label{lem:Topo2}
Consider a $(q+1)$-colored SYK graph $G$ of order $\delta$ with $S=\Psi_q(G)$, such that  for any $\C\in\llbracket 1, q\rrbracket$, $S_{\hat \C}$ is a collection of trees, and such that every core-chain contains an internal white chain-vertex. Then $G_{\hat 0}$ represents the connected sum of $\delta$   $(q-2)$-sphere bundles over $\cS^1$, $\cS^1\times\cS^{q-2}$.
\end{lemma}
\proof Every core-chain contains a white chain-vertex, which we denote by $v_{ij}$, $i$ and $j$ being the colors of the two incident edges that are not bridges. We pick a color $k\neq i,j$. From Corollary~\ref{cor:iResidue}, $S_{\hat k}$ is a forest, and since $S_{ij}$ is a sub-constellation of $S_{\hat k}$, it is a forest as well. From Lemma~\ref{lem:Forest-Bicolor}, $v_{ij}$ is therefore a handle-vertex. From Lemma~\ref{lem:Topo1}, $G_{\hat 0}$ represents a connected PL-manifold. Furthermore, the removal of $v_{ij}$ does not disconnect the graph: $S_{\setminus v_{ij}}$ is of excess $\delta-1$, so either it is a tree (in which case $(\Psi_q^{-1}(S_{\setminus v_{ij}}))_{\hat 0}$ represents a connected sphere from Prop.~\ref{MelonsAreTrees} and Prop.~\ref{prop:TopoMelo}), either $\delta-1>0$ and all of the remaining chains still contain a vertex of each color $\C\in\llbracket 1, q\rrbracket$, so that from Lemma~\ref{lemma:CSofSYK2}, $(\Psi_q^{-1}(S_{\setminus v_{ij}}))_{\hat 0}$ is still an SYK graph. We can apply Thm.~\ref{thm:Handles-Topo}, and proceed inductively on $S_{\setminus v_{ij}}$. \qed

\begin{theorem} 
\label{thm:Topo-OrderDelta}
Consider a large random bipartite $(q+1)$-edge-colored graph $G$ of fixed order~$\delta$. Then it is a.a.s.~dual to a unicellular space, and the triangulated boundary $\cT$ of the building block has the topology of the connected sum of $\delta$   $(q-2)$-sphere bundles over $\cS^1$, $\cS^1\times\cS^{q-2}$:
\be 
\lvert\cT\rvert \sim_{PL} \#_{\delta} (\cS^1\times\cS^{q-2}).
\ee 
\end{theorem}
\proof 
We have seen in Section~\ref{sec:SYK_conditions} that $G$ is a.a.s.~an SYK graph, which from Lemma~\ref{lemma:LongChainsColor} (Corollary~\ref{cor:iResidue}), a.a.s.~satisfies that for any $\C\in\llbracket 1, q\rrbracket$, $S_{\hat \C}$ is a collection of trees, with at least one internal white chain vertex per core-chain (Lemma~\ref{lemma:no_empty_white}). Therefore, a large $(q+1)$-edge-colored graph $G$ of fixed order~$\delta$ a.a.s.~satisfies the conditions (and thus the conclusions) of Lemma~\ref{lem:Topo2}. \qed

\

Note that because the residue $G_{\hat 0}$ almost surely does not represent a sphere, from Prop.~\ref{prop:PLMan}, the colored triangulation dual\footnote{It is obtained from the unicellular space described in Section~\ref{sub:Geom-Interp} by taking its cone.} to $G$ is almost surely not a piecewise-linear manifold, it has one singularity.

\section*{Concluding remarks}

As mentioned in the introduction, the colored tensor model and the SYK model share similar properties from a theoretical  physics perspective. On the level of the graphs, it had already been noticed in \cite{blt} that the perturbative expansions of the models in graphs differ, by comparison of the contributions to the first subleading orders. In the present work, we completed this analysis by performing a combinatorial study of the graphs at all orders.
This
 allows us to compare 
our results  
 with the 
 analysis performed by Gurau and Schaeffer in \cite{gilles} for the colored tensor model, i.e.~for bipartite $(q+1)$-edge-colored graphs of positive \emph{Gurau-degree}
\be 
\delta_\textrm{Gur}(G) = q + \frac {q(q-1)} 4 V(G) - \sum_{i<j =0}^q F_{ij}(G),
\ee 
where $F_{ij}(G)$ is the number of color-$ij$ cycles of the graph.

With the help of the bijection with constellations \cite{LL}, the analysis turns out to be considerably simpler for colored graphs of positive \emph{order} 
\begin{equation}
\deltaS_0(G) = 1+ \frac{q-1} 2 V(G) - F_0(G)
\end{equation}
than in the Gurau-Schaeffer case: in the case 
studied in this paper, 
we 
apply the 
method of kernel extraction. The main difficulty, 
as already noticed in \cite{blt, LL}, is to take into account the fact that SYK graphs are connected when deleting all color-0 edges. 
However, we prove here that
this turns out not to be a problem asymptotically: large random colored graphs of fixed order $\deltaS_0$ are a.a.s.~SYK graphs. 

While $(q+1)$-colored graphs of vanishing Gurau-degree (melonic graphs, Section~\ref{sub:Bij-Bip-case}) have a non-vanishing order $\deltaS_0=(q-1)(R_0 - 1)$, it is possible to show that the Gurau-degree of a large random $(q+1)$-edge-colored graph $G$ of fixed order~$\delta$ is a.a.s. 
\be 
\quad \delta_\mathrm{Gur}(G) =  q \delta 
\ee 
(and $\delta_\mathrm{Gur}(G_{\hat 0}) =  (q-1)\delta$ for its a.s.~connected 0-residue $G_{\hat 0}$). 
Indeed, at each step in the induction in the proof of Thm.~\ref{thm:Topo-OrderDelta}, the degree of $G$ (resp.~of $G_{\hat 0}$) decreases by $q$ (resp.~by $q-1$). These are illustrations of the differences between the two classifications of colored graphs: in terms of the Gurau-degree and in terms of the order. 

The differences between the two classifications,  as well as the simplicity of the present case, 
are better illustrated in the asymptotic enumerations of graphs of fixed Gurau-degree or of fixed order: while in the present case we obtain the 
estimate 
\begin{equation*}
\GF_{\delta}(z) =  P_\delta  \bigl( 1-\frac z {z_c}\bigr)^{\frac {1-3\delta}2}\Bigl[1 + O\Bigl(\sqrt{1-\frac z {z_c}}\Bigr)\Bigr], 
\end{equation*}
where the exponent $\frac {1-3\delta}2$ only depends on $\delta$,
the exponent in the asymptotic expression of the generating function of colored graphs of fixed Gurau-degree depends 
in a crucial way on $q$ (see  page 2 of \cite{gilles} with $q=D$), and it is only for $q>6$ that we obtain the analogous\footnote{In \cite{gilles}, for $q>6$, one has an analogous selection of ``trivalent'' schemes, while for $q<6$, the dominant schemes have a tree structure.} exponent $\frac{1-3\delta_\mathrm{Gur}/q}{2}$.

\section*{Acknowledgements}
 
The authors thank the two anonymous referees for their helpful comments and suggestions, and Valentin Bonzom for interesting discussions. 
Adrian Tanasa is partially supported by the grants CNRS Infiniti "ModTens" and PN 09 37 01 02. Luca Lionni is a JSPS international research fellow.
\'Eric Fusy is partially supported by the CNRS Infiniti "ModTens" grant and by the ANR grant GATO, ANR-16-CE40-0009.

\bigskip

\noindent \'Eric Fusy\\
{\it\small  LIX, CNRS UMR 7161, \'Ecole Polytechnique, 91120 Palaiseau, France, EU.}

\medskip

\noindent Luca Lionni\\
{\it\small  Yukawa Institute for Theoretical Physics, Kyoto University, Japan.}

\medskip

\noindent
Adrian Tanasa\\
{\it\small LABRI, CNRS UMR 5800, 
Universit\'e Bordeaux,} {\it\small 351 cours de la Lib\'eration, 33405 Talence cedex, France, EU.}\\
{\it\small Horia Hulubei National Institute for Physics and Nuclear Engineering,
P.O.B. MG-6, 077125 Magurele, Romania, EU.}\\
{\it\small I. U. F., 1 rue Descartes, 75005 Paris, France, EU.}
\end{document}